\documentclass{article}
\usepackage{amsmath,amssymb}

\newcommand{\dis}{\displaystyle}

\newcommand{\DD}{\mathcal{D}}
\newcommand{\FF}{\mathcal{F}}
\newcommand{\Q}{\overline{Q}}

\newcommand{\Bu}{\overline{B}}
\newcommand{\Bb}{\underline{B}}

\newcommand{\pn}{\overline{\psi}_n}

\newcommand{\B}{\overline{B}}
\newcommand{\Div}{\mathrm{div}}

\newcommand{\R}{\mathbb{R}}

\newcommand{\eps}{\varepsilon}
\newcommand{\dt}{\partial_t}

\newcommand{\ve}{v_\eps}
\newcommand{\Pe}{P_\eps}

\def\qed{\hbox{${\vcenter{\vbox{
  \hrule height 0.4pt\hbox{\vrule width 0.4pt height 6pt
  \kern5pt\vrule width 0.4pt}\hrule height 0.4pt}}}$}}

\newtheorem{theo}{Theorem}
\newtheorem{prop}[theo]{Proposition}
\newtheorem{lemm}[theo]{Lemma}

\date{}
\author{Alexis Vasseur \thanks{Department of Mathematics, University of Texas}}
\title{Higher derivatives estimate for the 3D Navier-Stokes equation}
\begin{document}
\maketitle
\bibliographystyle{plain}
\noindent{\bf Abstract:} In this article, a non linear family of spaces, based on the energy dissipation, is introduced.
This family bridges an  energy space  (containing weak solutions to Navier-Stokes equation) to a critical space (invariant through the canonical scaling of the Navier-Stokes equation). This family is used to get  uniform estimates on higher derivatives  to  solutions to the 3D Navier-Stokes equations. Those estimates are uniform, up to the possible blowing-up time. The proof uses
blow-up techniques. Estimates can be obtained by this means thanks to the galilean
invariance of the transport part of the equation.



\section{Introduction}

In this paper, we investigate estimates of higher derivatives of
solutions to the incompressible Navier-Stokes equations in dimension
3, namely:
\begin{equation}\label{eq_NS}
\begin{array}{l}
\dis{\dt u+\Div(u\otimes u)+\nabla P-\Delta u=0\qquad
t\in (0,\infty), \ x\in\R^3,}\\[3mm]
\dis{\mathrm{div} u=0.}
\end{array}
\end{equation}
The initial  value problem is endowed with the conditions:
\begin{eqnarray*}
&&u(0,\cdot)=u^0\in L^2(\R^3).
\end{eqnarray*}

The existence of weak solutions for this problem was proved long ago
by Leray \cite {Leray} and Hopf \cite{Hopf}. For this, Leray
introduces a notion of weak solution. He shows that for any initial
value with finite energy $u^0\in L^2(\R^3)$ there exists a function
$u\in L^\infty (0,\infty;L^2(\R^3))\cap
L^2(0,\infty;\dot H^1(\R^3))$ verifying (\ref{eq_NS}) in the sense of
distribution. From that time on, much effort has been made to
establish results on the uniqueness and regularity of weak
solutions. However those two questions remain yet mostly open.
Especially it is not known until now if such a weak solution can
develop singularities in finite time, even considering smooth
initial data. We present our main result on a laps of time $(0,T)$ where the solution
is indeed smooth (with possible blow-ups both at $t=0$ and $t=T$). We will carefully show, however, that the estimates do not depend on the blow-up time $T$, but only on
$\|u^0\|_{L^2}$ and $\inf(t,1)$.
 The aim of
this paper is to show the following theorem.
\begin{theo}\label{theo_main}
For  any $t_0>0$, any $\Omega$ bounded subset of $(t_0,\infty)\times \R^3$, any integer $n\geq1$, any $\gamma>0$,
and any $p\geq0$ such that
\begin{equation}\label{scale_nonlinear}
\frac{4}{p}>n+1,
\end{equation}
there exists a constant $C$, such that the following property holds.

For any smooth solution $u$ of (\ref{eq_NS}) on $(0,T)$ (with possible blow-up at 0 and $T$), we have
$$
\|\nabla^n u\|_{L^p(\Omega\cap[(0,T)\times\R^3])}\leq C\left(\|u^0\|^{2(1+\gamma)/p}_{L^2(\R^3)}+1\right).
$$
Note that the constant $C$ does not depend on the solution $u$ nor on the blowing-up time $T$.
\end{theo}
Note that for $n\geq 3$ we consider $L^p$ spaces with $p<1$. Those spaces are not complete. For this reason the result cannot be easily extend to general weak solutions after the possible blow-up time. However, up to $d= 2$, the result can be proven in this context.
For this reason, along the proof, we will always consider suitable weak solutions, following \cite{CKN}.  That is, solutions verifying in addition to
(\ref{eq_NS}) the generalized energy inequality in the sense of
distribution:
\begin{equation}\label{eq_NSenergy}
\dt \frac{|u|^2}{2}+\Div\left(u\frac{|u|^2}{2}\right)+\Div (uP)+|\nabla
u|^2-\Delta \frac{|u|^2}{2}\leq0\qquad t\in (0,\infty), \
x\in\R^3.
\end{equation}
Moreover, by interpolation, 
the result of Theorem \ref{theo_main} can be extended to the whole real derivative coefficients, $1<d\leq 2$, for $\|\Delta^{d/2}u\|_{L^{p}}$ with
$$
\frac{4}{p}>d+1.
$$

Our result can be seen as a
kind of anti-Sobolev result. Indeed, as we will see later,
 $\|\nabla u\|^2_{L^2}$ is used as a pivot quantity to  control higher derivatives on the solution. The result for $d=2$ was
 obtained in a slightly better space, with completely different techniques by Lions
 \cite{Lions}.  He shows that $\nabla^2 u$ can be bounded in the Lorentz space $L^{4/3,\infty}$.

In a standard way, using the energy inequality and interpolation, we get estimates on
$\Delta^{d/2}u\in L^p((0,\infty)\times\R^3)$ for
\begin{equation}\label{scale_energy}
\frac{5}{p}=d+\frac{3}{2},\qquad 0\leq d\leq1.
\end{equation}
The Serrin-Prodi conditions (see \cite{Serrin},\cite{Fabes},
\cite{Struwe}) ensure the regularity for solutions such that
$\Delta^{d/2}u\in L^p((0,\infty)\times\R^3)$ for
\begin{equation}\label{scale_regularity}
\frac{5}{p}=d+1,\qquad 0\leq d<\infty.
\end{equation}
Those two families of spaces  are given by an affine relation on $d$  with
respect to $1/p$ with slope $5$. Notice that the family of spaces present in
Theorem \ref{theo_main} has a different slope. Imagine, that we were able to extend this result along the same line
with $d<1$. For $d=0$, we would obtain almost $u\in
L^4((0,\infty)\times\R^3)$, which would imply that the energy
inequality (\ref{eq_NSenergy}) is an equality (see \cite{Temam}).
Notice also that the line of this new family of spaces crosses the line of the
critical spaces (\ref{scale_regularity}) at $d=-1$, $1/p=0$. This
point corresponds (at least formally) to the Tataru and Koch result
on regularity of solutions small in
$L^\infty(0,\infty;BMO^{-1}(\R^3))$ (see \cite{Tataru}). However,
at this time, due to the ``anti-Sobolev" feature of the proof,
obtaining results for $d<1$ seems out of reach.

To see where lie the difficulties, let us focus on the result on the
third derivatives. Consider the gradient of the Navier-Stokes
equations (\ref{eq_NS}).
$$
\dt \nabla u-\Delta\nabla u=  -\nabla u\cdot\nabla
u-\nabla^2P-(u\cdot\nabla)\nabla u.
$$
Note that the two first right-hand side terms lie  in
$L^1((0,\infty)\times\R^3)$ (for the pressure term, see
\cite{Lions}). Parabolic regularity are not complete in $L^1$. This
justify the fact that we miss the limit case  $L^{1}$. But, surprisingly, the worst term
is the transport one $(u\cdot\nabla)\nabla u$. To control it in
$L^1$ using the control on $D^2u$ in $L^{4/3,\infty}$ of Lions
\cite{Lions}, we would need $u\in L^{4,1}$, which is not known.
To overcome this difficulty, we will consider the solution in another frame, locally,
by following the flow.

The idea of the proof comes from the result of partial regularity
obtained by Caffarelli, Kohn and Nirenberg \cite{CKN}. This paper
extended  the analysis about the possible singular points set,
initialized by Scheffer in a series of paper
\cite{Scheffer1,Scheffer2,Scheffer3,Scheffer4}. The main remark in
\cite{CKN} is that the dissipation of entropy
\begin{equation}\label{eq_diss_entropy}
\DD(u)=\int_0^\infty\int_{\R^3}|\nabla u|^2\,dx\,dt
\end{equation}
has a scaling, through the standard invariance of the equation,
which is far more powerful that any other quantities from the energy
scale (\ref{scale_energy}). Let us be more specific. The standard
invariance of the equation gives that for any $(t_0,x_0)\in
\R^+\times\R^3$ and $\eps>0$, if $u$ is a suitable solution of the
Navier-Stokes equations (\ref{eq_NS}) (\ref{eq_NSenergy}), then
\begin{equation}\label{eq_invariance}
u_\eps(t,x)=\eps u(t_0+\eps^2 t,x_0+\eps x)
\end{equation}
is also solution to (\ref{eq_NS}) (\ref{eq_NSenergy}). The
dissipation of energy of this quantity is then given by
$$
\DD(u_\eps)=\eps^{-1}\DD(u).
$$
This power of $\eps$ made possible in \cite{CKN} to show that the
Hausdorff dimension of the set of blow-up points is at most 1. This
was a great improvement of the result obtained by Scheffer who gives
5/3 as an upper bound for the Hausdorff dimension of this set. We can notice
that it is what we get considering the quantity of the energy scale
(\ref{scale_energy}) with $d=0, p=10/3$:
$$
\FF(u)=\int_0^\infty\int_{\R^3}|u|^{10/3}\,dx\,dt.
$$
Indeed:
$$
\FF(u_\eps)=\eps^{-5/3}\FF(u).
$$
The idea of this paper is to give a quantitative version of the
result of \cite{CKN}, in the sense, of getting control of norms of
the solution which have the same nonlinear scaling that $\DD$.
Indeed, for any norm of the non linear scaling
(\ref{scale_nonlinear}), we have (in the limit case)
$$
\|\nabla ^nu_\eps\|^p_{L^{p}}=\eps^{-1}\|\nabla^nu\|^p_{L^{p}}.
$$
The paper is organized as follows. In the next section, we give some preliminaries and fix some notations. We introduce the local frame following the flow in the third section. The fourth section is dedicated to a local result providing a universal control of the higher derivatives  of $u$  from a local control of the dissipation of the energy $\|\nabla u\|^2_{L^2}$ and a corresponding quantity on the pressure (see Proposition \ref{prop_estimates}). Ideally, we would like to consider a quantity on the pressure which has
the same nonlinear scaling as $\DD(u)$. The corresponding quantity is
$\|\nabla^2 P\|_{L^1}$. Unfortunately, we need a slightly better integrability in time for the local study. This is the reason why we miss the limit case $L^{p,\infty}$ with
$$
\frac{4}{p}=n+1.
$$
This is also the reason why we need to work with fractional Laplacian for the pressure:
$\|\Delta^{-s} \nabla^2P\|_{L^p}$ with $0<s<1/2$.
In the last section, we show how this local study leads to our main theorem.

\section{Preliminaries and notations}

Let us denote $Q_r=(-r^2,0)\times B_r$ where $B_r=B(0,r)$, the ball in
$\R^3$ of radius $r$ and centered at 0.

For $F\in L^p(\R^+\times\R^3)$, we define the Maximal function in $x$ only by
$$
MF(t,x)=\sup_{r>0}\frac{1}{r^3}\int_{B_r}|F(t,x+y)|\,dy.
$$
We recall that for any $1<p<\infty$, there exists $C_p$ such that for any $F\in L^p(\R^+\times\R^3)$
$$
\|MF\|_{L^p(\R^+\times\R^3)}\leq C_p \|F\|_{L^p(\R^+\times\R^3)}.
$$
Moreover, there exists a constant $C$ such that  for any  $F\in L^1(\R^+;\mathcal{H}(\R^3))$, (where $\mathcal{H}$ stands for the Hardy space), then
$$
\|MF\|_{L^1(\R^+\times\R^3)}\leq C \|F\|_{L^1(\R^+;\mathcal{H}(\R^3))}.
$$
We begin with an interpolation lemma. It is a straightforward consequence of a result in \cite{Bergh}. We state it here for further reference.

\begin{lemm}\label{lemm_interpolation}
For any function F such that $(-\Delta)^{d_1/2}F$ lies in $L^{p_1}(0,\infty;L^{q_1}(\R^3))$
and $(-\Delta)^{d_2/2}F\in L^{p_2}(0,\infty;L^{q_2}(\R^3))$ with
$$
d_1,d_2\in\R,\qquad 1\leq p_1,p_2\leq\infty, \qquad 1< q_1,q_2<\infty,
 $$
we have $(-\Delta)^{d/2}F\in L^{p}(0,\infty;L^q(\R^3))$
with
\begin{eqnarray*}
&&\|(-\Delta)^{d/2}F\|_{L^{p}(0,\infty;L^q(\R^3))}\\
&&\qquad\leq \|(-\Delta)^{d_1/2}F\|^\theta_{L^{p_1}(0,\infty;L^{q_1}(\R^3))}\|(-\Delta)^{d_2/2}F\|^{1-\theta}_{L^{p_2}(0,\infty;L^{q_2}(\R^3))},
\end{eqnarray*}
for any $d,p,q$ such that
\begin{eqnarray*}
&&\frac{1}{q}=\frac{\theta}{q_1}+\frac{1-\theta}{q_2},\\
&&\frac{1}{p}=\frac{\theta}{p_1}+\frac{1-\theta}{p_2},\\
&&d=\theta d_1+(1-\theta) d_2,
\end{eqnarray*}
where $0< \theta< 1$.
\end{lemm}

\noindent{\bf Proof.}
Exercise 31 page 168 in \cite{Bergh} shows that for any $0< t< \infty$, we have
$$
\|(-\Delta)^{d/2}F(t)\|_{L^{p}(\R^3)}\leq \|(-\Delta)^{d_1/2}F(t)\|^\theta_{L^{p_1}(\R^3)}\|(-\Delta)^{d_2/2}F(t)\|^{1-\theta}_{L^{p_2}(\R^3)}.
$$
Interpolation in the time variable gives the result.
\qquad\qed
\vskip0.3cm
In the second lemma we show that we can control a local $L^1$ norm on a function $f$  by its mean value and some local control on the maximal function of $(-\Delta)^{-s}\nabla f$, $0<s<1/2$.
This extends the fact that we can control the local $L^1$ norm by the mean value and a local $L^p$ norm of the gradient. But due to the nonlocal feature of the fractional Laplacian, we need to consider the maximal function to recapture all the information needed.

\begin{lemm}\label{lemm_local gradient control}
Let $0<s<1/2$, $q\geq 1$, $p\geq 1$. For any $\phi\in C^\infty(\R^3)$, $\phi\geq0$, compactly supported in $B_{1}$ with $\int_{\R^3}\phi(x)\,dx=1$, there exists $C>0$ such that, for any function $f\in L^q(\R^3)$ with $(-\Delta)^{-s}\nabla f\in L^p(\R^3)$ and $|\int f\phi\,dx|$ bounded, we have $f\in L^1(B_1)$ and
$$
\|f\|_{L^1(B_1)}\leq C\left(\left|\int_{\R^3} f(x)\phi(x)\,dx\right|+\|M((-\Delta)^{-s}\nabla f)\|_{L^p(B_1)}\right).
$$
\end{lemm}
\noindent{\bf Proof.}
Let us denote $g=(-\Delta)^{-s}\nabla f$.
Since $f\in L^q(\R^3)$, we have
$$
f=-(-\Delta)^{s-1}\Div g.
$$
So, for any $x\in B_1$
$$
f(x)=C_s \int_{\R^3}\frac{g(y)}{|x-y|^{2(1+s)}}\cdot\frac{(x-y)}{|x-y|}\,dy,
$$
and
\begin{eqnarray*}
&&f(x)-\int_{\R^3}\phi(z)\,f(z)\,dz\\
&&\qquad=C_s \int_{\R^3}\int_{\R^3}\phi(z)g(y)\left(\frac{(x-y)/|x-y|}{|x-y|^{2(1+s)}}-\frac{(z-y)/|z-y|}{|y-z|^{2(1+s)}}\right)\,dy\,dz.
\end{eqnarray*}

Note that,  for $k\geq 2$, $y\in B_{2^k}\setminus B_{2^{k-1}}$, $x\in B_1$, $z\in B_1$, we have
$$
\left|\frac{(x-y)/|x-y|}{|x-y|^{2(1+s)}}-\frac{(z-y)/|z-y|}{|y-z|^{2(1+s)}}\right|\leq \frac{C}{2^{k(3+2s)}}.
$$
Moreover
\begin{eqnarray*}
&&\int_{B_1}\int_{B_1}\int_{B_2}\phi(z)|g(y)|\left|\frac{(x-y)/|x-y|}{|x-y|^{2(1+s)}}-\frac{(z-y)/|z-y|}{|y-z|^{2(1+s)}}\right|\,dy\,dz\,dx\\
&&\qquad \leq\int_{B_3}\int_{B_1}\int_{B_2}\frac{\phi(z)|g(y)|}{|x|^{2(1+s)}}\,dy\,dz\,dx+\int_{B_1}\int_{B_3}\int_{B_2}\frac{\sup|\phi||g(y)|}{|z|^{2(1+s)}}\,dy\,dz\,dx\\
&&\qquad \leq 2C_s \|g\|_{L^1(B_1)}\leq 2C_s \|Mg\|_{L^1(B_1)},
\end{eqnarray*}
since $2(1+s)<3$.
Hence
\begin{eqnarray*}
&&\qquad\qquad \left\|f-\int\phi(z)f(z)\,dz\right\|_{L^1(B_1)}\\
&&\leq \int_{B_1}\int_{B_1}\int_{B_2}\phi(z)|g(y)|\left|\frac{(x-y)/|x-y|}{|x-y|^{2(1+s)}}-\frac{(z-y)/|z-y|}{|y-z|^{2(1+s)}}\right|\,dy\,dz\,dx\\
&&\qquad +\sum_{k=2}^\infty\int_{B_1}\int_{B_1}\int_{(B_{2^k}\setminus B_{2^{k-1}})}\!\!\!\!\!\!\!\!\!\!\!\phi(z)|g(y)|\left|\frac{(x-y)/|x-y|}{|x-y|^{2(1+s)}}-\frac{(z-y)/|z-y|}{|y-z|^{2(1+s)}}\right|\\
&&\leq 2C_s\|Mg\|_{L^1(B_1)}+C\sum_{k=2}^\infty \int_{B_{2^k}}\frac{|g(y)|}{2^{k(3+2s)}}\,dy\\
&&\leq 2C_s\|Mg\|_{L^1(B_1)}+8C\sum_{k=2}^\infty 2^{-2s k}\frac{1}{|B_{2^{k+1}}|}\int_{B_1}\int_{B_{2^{k+1}}}|g(y+u)|\,dy\,du\\
&&\leq 2C_s\|Mg\|_{L^1(B_1)}+C\|Mg\|_{L^1(B_1)}\sum_{k=2}^\infty[2^{-2s}]^k\\
&&\leq C_s\|Mg\|_{L^1(B_1)},
\end{eqnarray*}
whenever $0<s<1/2$.\qquad \qed

\vskip0.3cm
We give now  very standard results of parabolic regularity. There are not even optimal, but enough for our study.
\begin{lemm}\label{lemm_parabolic regularity}
For any $1<p<\infty$, $t_0>0$, there exists a constant $C$ such that the following is true.
Let $f,g\in L^p((-t_0,0)\times \R^3)$ be compactly supported in $B_1$. Then there exists a unique $u\in L^p(-t_0,0;W^{1,p}(\R^3))$  solution  to
\begin{eqnarray*}
&&\dt u-\Delta u=g+\Div f, \qquad -t_0\leq t\leq 0,\ \  x\in \R^3,\\
&& u(-t_0,x)=0,\qquad x\in \R^3.
\end{eqnarray*}
Moreover, 
\begin{equation}\label{eq_4}
\|u\|_{L^p(-t_0,0;W^{1,p}(B_1))}\leq C (\|f\|_{L^p((-t_0,0)\times\R^3)}+\|g\|_{L^p((-t_0,0)\times\R^3)}).
\end{equation}
If $g\in L^1(-t_0,0;L^\infty(\R^3))$ and $f\in L^1(-t_0,0;W^{1,\infty}(\R^3))$,
then
$$
\|u\|_{L^\infty(-t_0,0)\times\R^3)}\leq C(\|g\|_{L^1(-t_0,0;L^\infty(\R^3))}+\|f\|_{L^1(-t_0,0;W^{1,\infty}(\R^3))}).
$$
\end{lemm}

\noindent{\bf Proof.}
We get the solution using the Green function:
$$
u(t,x)=\int_{-t_0}^t\frac{1}{4\pi(t-s)^{3/2}}\int_{\R^3}e^{-\frac{|x-y|^2}{4(t-s)}}(g(s,y)+\Div f(s,y))\,dy\,ds.
$$
From this formulation, using that $z^n e^{-z^2}$ are bounded functions,  we find that
\begin{equation}\label{eq_6}
|u(t,x)|\leq C\frac{\|f\|_{L^1((-t_0,0)\times B_1)}+\|g\|_{L^1((-t_0,0)\times B_1)}}{|x|^3},\qquad \mathrm{for} \ |x|>2, -t_0\leq t<0.
\end{equation}
Standard  Solonnikov's parabolic regularization result  gives (\ref{eq_4}) (see for instance \cite{Solonnikov}).
Finally, if  $g\in L^1(-t_0,0;L^\infty(\R^3))$ and $f\in L^1(-t_0,0;W^{1,\infty}(\R^3))$, then the function
$$
v(t,x)=\int_0^t(\|g(s)\|_{L^\infty}+\|\Div f(s)\|_{L^\infty})\,ds
$$
is a supersolution thanks to (\ref{eq_6}). The global bound follows.
\qquad\qed

\vskip0.3cm
The last lemma of this section  is a standard  decomposition of  the pressure term as a close range part and a long range part.
\begin{lemm}\label{lemm_decomposition of pressure}
Let $\Bu$ and $\Bb$ be two balls such that
$$
\Bu\subset\subset \Bb.
$$
Then for any $1<p<\infty$, there exists a constant $C>0$ and a family of constants $\{C_{d,q}\ \setminus \ d,q\ \ \mathrm{integers}\}$ (depending only on $p$, $\Bb$ and $\Bu$) such that for any $R\in L^1(\Bb)$ and $A\in [L^p(\Bb)]^{N\times N}$ symmetric matrix, verifying
$$
-\Delta R=\Div\Div A,\qquad\mathrm{in}\ \ \Bb,
$$
we have a decomposition
$$
R=R_1+R_2,
$$
with, for any  integer $q\geq0$, $d\geq0$:
\begin{eqnarray*}
&&\|R_1\|_{L^p(\Bu)}\leq C \|A\|_{L^p(\Bb)},\\
&&\|\nabla^d R_2\|_{L^\infty(\Bu)}\leq C_{d,q}\left(\|A\|_{L^1(\Bb)}+\|R\|_{W^{-q,1}(\Bb)}\right).
\end{eqnarray*}
Moreover, if $A$ is Lipschitzian, then we can choose $R_1$ such that
$$
\|R_1\|_{L^\infty(\Bu)}\leq C\left(\|\nabla A\|_{L^\infty(\Bb)}+\| A\|_{L^\infty(\Bb)}\right).
$$
\end{lemm}
\noindent{\bf Proof.}

Let $B^*$ be a a ball such that
$$
\Bu\subset\subset B^*\subset\subset\Bb,
$$
with a distance between $\Bu$ and ${B^*}^c$ bigger that $D/2$, where $D$ is the distance between $\Bu$ and $\Bb^c$.
Consider a smooth nonnegative  cut-off function $\psi$, $0\leq \psi\leq 1$ such that
\begin{eqnarray*}
\psi(x)&=&1\qquad \mathrm{in} \ B^*,\\
&=&0\qquad \mathrm{in} \ \Bb^c.
\end{eqnarray*}
Then the function $\psi R$ (defined in $\R^3$) is solution in $\R^3$ to
\begin{eqnarray*}
-\Delta (\psi R)&=&\Div\Div (\psi A)\\
&&+R\Delta \psi+A:\nabla^2\psi\\
&&-2\Div\{\nabla\psi\cdot A+R\nabla\psi\}.
\end{eqnarray*}
We denote
\begin{eqnarray*}
&&R_1=(-\Delta)^{-1}\Div\Div(\psi A),\\
&&R_2=(-\Delta)^{-1}\left(R\Delta \psi+A:\nabla^2\psi-2\Div\{\nabla\psi\cdot A+R\nabla\psi\}\right).
\end{eqnarray*}
We have, on $\Bu$, $R=R_1+R_2$. The operator $(-\Delta)^{-1}\Div\Div$ is a Riesz operator, so there exists a constant (depending only on $p$ and $\psi$) such that
\begin{eqnarray*}
&&\|R_1\|_{L^p(\R^3)}\leq C \|\psi A\|_{L^p(\R^3)}\leq C\|A\|_{L^p(\Bb)},\\
&&\|R_1\|_{C^\alpha(\R^3)}\leq C \|\psi A\|_{C^\alpha(\R^3)}\leq C\left(\|\nabla A\|_{L^\infty(\Bb)}+\|A\|_{L^\infty(\Bb)}\right).
\end{eqnarray*}
Using the fact that $\nabla \psi$ and $\nabla^2\psi$ vanishes on $B^*\cup \Bb^c$, we have for any $x\in \Bu$:
\begin{eqnarray*}
&&|\nabla^d R_2(x)|
=\left|\int_{\R^3}\nabla^d\left(\frac{1}{|x-y|}\right)\left(R\Delta \psi+A:\nabla^2\psi\right)(y)\,dy\right.\\
&&\qquad\qquad\qquad\left.  +2\int_{\R^3}\nabla^{d+1}\left(\frac{1}{|x-y|}\right)\{\nabla\psi\cdot A+R\nabla\psi\}(y)\,dy\right|\\
&&\qquad\leq \|\nabla^2\psi\|_{L^\infty}\|A\|_{L^1(\Bb)}\sup_{|x-y|\geq D/2}\left|\nabla^d\left(\frac{1}{|x-y|}\right)\right|\\
&&\qquad\qquad +2\|\nabla\psi\|_{L^\infty}\|A\|_{L^1(\Bb)}\sup_{|x-y|\geq D/2}\left|\nabla^{d+1}\left(\frac{1}{|x-y|}\right)\right|\\
&&\qquad\qquad +\|R\|_{W^{-q,1}(\Bb)}\sup_{|x-y|\geq D/2}\left|\nabla^q\left[\nabla^d\left(\frac{1}{|x-y|}\right)\Delta \psi\right]\right|\\
&&\qquad\qquad +2\|R\|_{W^{-q,1}(\Bb)}\sup_{|x-y|\geq D/2}\left|\nabla^q\left[\nabla^{d+1}\left(\frac{1}{|x-y|}\right)\nabla \psi\right]\right|\\
&&\qquad\leq C_{d} \left[\left(\frac{2}{D}\right)^{d+2}+\left(\frac{2}{D}\right)^{d+1}\right]\|A\|_{L^1(\Bb)}\\
&&\qquad\qquad +C_{d,q} \left[\left(\frac{2}{D}\right)^{d+1}+\left(\frac{2}{D}\right)^{q+d+2}\right]\|R\|_{W^{-q,1}(\Bb)}.
\end{eqnarray*}
\qquad \qed

\section{Blow-up method along the trajectories}

Our result relies on a  local study, which was the keystone of the partial regularity result of \cite{CKN}. (see \cite{Lin} for an other proof). We use, here,  the version of  \cite{Vasseur}. This version is better for our purpose because it requires a bound on the pressure only in $L^p$ in time for any $p>1$.
\begin{prop}\label{prop_main}\cite{Vasseur}
 For any $p>1$,
there exists $\eta>0$,   such that the following property holds. For any
$u$, suitable weak solution to the Navier-Stokes equation
(\ref{eq_NS}),  (\ref{eq_NSenergy}), in $Q_1$,
such that
\begin{eqnarray}
&&\sup_{-1<t<0}\!\!\left(\int_{B_1}\!\!|u(t,x)|^2\,dx\right)+\!\int_{Q_1}\!\!|\nabla u|^2\,dx\,dt\!+\!\int_{-1}^0\!\!\left(\int_{B_1}\!\!|P|\,dx\right)^p \!dt\leq
\eta,\label{Hypenergiepetite}
\end{eqnarray}
we have
$$
\sup_{(t,x)\in Q_{1/2}}|u(t,x)|\leq 1.
$$
\end{prop}
As explained in the introduction, the proof of Theorem \ref{theo_main} relies on this local control. From there  we can get control on  higher derivatives of $u$.
We first show the following lemma. It introduces the pivot quantity. Note that the ideal pivot quantity would be $\|\nabla u\|^2_{L^2(L^2)}+\|\nabla^2 P\|_{L^1(L^1)}$. This is because this quantity scales as $1/\eps$ through the canonical scaling. However, to use Proposition \ref{prop_main} locally, we need a better integrability in time on the pressure. For this reason, we add the quantity on the pressure involving the fractional Laplacian. We get a better integrability in time on the pressure, at the cost of a slightly worst
rate of change in $\eps$ through the canonical scaling. Finally, due to the nonlocal character of the fractional Laplacian, the maximal function is used in order to recapture all the local information needed (see Lemma \ref{lemm_local gradient control}).

\begin{lemm}\label{lemmmaximalenergy}
For   any $0<\delta<1$, there exists $\gamma>0$ and  a constant $C>0$ such that for any $u$ solution to (\ref{eq_NS}) (\ref{eq_NSenergy}), with $u^0\in L^2(\R^3)$, we have
\begin{eqnarray*}
&&\int_0^\infty\int_{\R^3}\left(|M((-\Delta)^{-\delta/2}\nabla^2 P)|^{1+\gamma}+|\nabla^2 P|+|\nabla u|^{2}\right)\,dx\,dt\\
&&\qquad\qquad\leq C\left(\|u^0\|^2_{L^2(\R^3)}+\|u^0\|^{2(1+\gamma)}_{L^2(\R^3)}\right).
\end{eqnarray*}
Moreover, $\gamma$ converges to 0 when $\delta$ converges to 0.
\end{lemm}
\noindent{\bf Proof.}
Integrating in $x$ the energy equation (\ref{eq_NSenergy}) gives that
 \begin{equation}\label{eq_2}
 \int_0^\infty\int_{\R^3}|\nabla u|^{2}\,dx\,dt\leq \|u^0\|^2_{L^2(\R^3)},
 \end{equation}
 together with
 $$
 \|u\|^2_{L^\infty(0,\infty;L^2(\R^3))}\leq \|u^0\|^2_{L^2(\R^3)}.
 $$
 By Sobolev imbedding and interpolation, this gives in particular
 that
 \begin{equation}\label{eq_1}
 \|u\|^2_{L^4(0,\infty;L^3(\R^3))}\leq  C\|u^0\|^2_{L^2(\R^3)}.
 \end{equation}
 For the pressure, we have $\nabla^2 P\in L^1(\mathcal{H})$ (see Lions \cite{Lions}).
Indeed,
\begin{eqnarray*}
&&\nabla^2 P=(\nabla^2\Delta^{-1})\sum_{ij}\partial_i u_j\partial_ju_i\\
&&\qquad=(\nabla^2\Delta^{-1})\sum_i(\partial_i u)\cdot\nabla u_i.
\end{eqnarray*}
For any $i$, we have $\mathrm{rot}(\nabla u_i)=0$ and $\Div\ \partial_i u=0$. Hence, from the div-rot lemma (see Coifman, Lions, Meyer and Semmes \cite{Coifman}), we have
$$
\|\sum_i \partial_i u\cdot\nabla u_i\|_{L^1(\mathcal{H})}\leq \|\nabla u\|^2_{L^2}.
$$
But $\nabla^2\Delta^{-1}$ is a Riesz operator (in $x$ only) which is bounded from $\mathcal{H}$ to $\mathcal{H}$.
Hence:
\begin{equation}\label{eq_bonne pression}
\|\nabla^2 P\|_{L^1(\R^+\times\R^3)}\leq C\|\nabla^2 P\|_{L^1(\R^+;\mathcal{H}(\R^3))}\leq C\|\nabla u\|^2_{L^2(\R^+\times\R^3)}.
\end{equation}
By Sobolev imbedding, for any $0<s<1$, we have
\begin{equation}\label{eq_5}
\|(-\Delta)^{-s/2}\nabla^2 P\|_{L^{1}(0,\infty;L^{p}(\R^3))}\leq C \|u^0\|^2_{L^2}
\end{equation}
for
$$
\frac{1}{p}=1-\frac{s}{3}.
$$
we have also
$$
(-\Delta)^{-1/2}\nabla^2 P=\sum_{ij}[(-\Delta)^{-3/2}\nabla^2\partial_i](\partial_j u_i u_j).
$$
The operators $(-\Delta)^{-3/2}\nabla^2\partial_i$ are Riesz operators so, together with (\ref{eq_2}) (\ref{eq_1}), we have
\begin{equation}\label{eq_3}
\|(-\Delta)^{-1/2}\nabla^2 P\|_{L^{4/3}(0,\infty;L^{6/5}(\R^3))}\leq C \|u^0\|^2_{L^2(\R^3)}.
\end{equation}
By interpolation with (\ref{eq_5}), using  Lemma \ref{lemm_interpolation} with $\theta=1/(1+4s)$, we find
$$
\|M[(-\Delta)^{-\delta/2}\nabla^2 P]\|_{L^{1+\gamma}((0,\infty)\times\R^3)}\leq C\|u^0\|^2_{L^2(\R^3)}
$$
with
$$
\delta=\frac{5s}{1+4s},\qquad\qquad
 \gamma=\frac{s}{1+3s}.
$$

Note that $\gamma$ converges to 0 when $\delta$ goes to 0. This, together with (\ref{eq_bonne pression}) and (\ref{eq_2}), gives the result.
\qquad \qed

\vskip0.3cm

Let us fix from now on a smooth cut-off function $0\leq\phi\leq 1$ compactly supported in $B_1$ and such that
\begin{equation}\label{etoile1}
\int_{\R^3}\phi(x)\,dx=1.
\end{equation}
 For any $\eps>0$, we define
\begin{equation}\label{etoile2}
u_\eps(t,x)=\int_{\R^3}\phi(y)u(t,x+\eps y)\,dy.
\end{equation}
Note that $u_\eps\in L^\infty(0,\infty;C^\infty(\R^3))$ and $\Div
u_\eps=0$. We define the flow:
\begin{equation}\label{eq_flow}
\begin{array}{l}
\dis{\frac{\partial X}{\partial s}=u_\eps(s,X(s,t,x))}\\[0,3cm]
\dis{X(t,t,x)=x.}
\end{array}
\end{equation}
Consider, for any $0<\delta<1$ and $\eta^*>0$:
$$
\Omega^\delta_\eps=\left\{(t,x)\in (4\eps^2,\infty)\times\R^3\ | \
\frac{1}{\eps}\int_{t-4\eps^2}^{t}\int_{B_{2\eps}}\!\!\!\!F^\delta(s,X(s,t,x)+y)\,ds\,dy\leq \eta^*\eps^{\delta}\right\},
$$
where
$$
F^\delta(t,x)=|M((-\Delta)^{-\delta/2}\nabla^2 P)|^{1+\gamma}+|\nabla u|^{2}+|\nabla^2 P|,
$$
and $\gamma$ is defined in Lemma \ref{lemmmaximalenergy}.
We then have the following lemma.
\begin{lemm}\label{lemm_main}
There exists a constant $C$ such that for any $0<\eps<1$, $0<\delta<1$, and $\eta^*>0$ we have
$$
|[\Omega^\delta_\eps]^c|\leq C\left(\frac{\|u^0\|^2_{L^2(\R^3)}+\|u^0\|^{2(1+\gamma)}_{L^2(\R^3)}}{\eta^*}\right)\eps^{4-\delta}.
$$
\end{lemm}
\noindent{\bf Proof.} Define for $t>4\eps^2$
\begin{equation}\label{eq_F}
F^\delta_\eps(t,x)=\frac{1}{(2\eps)^5}\int_{t-4\eps^2}^{t}\int_{B_{2\eps}}F^\delta(s,X(s,t,x)+y)\,ds\,dy.
\end{equation}
We have
\begin{eqnarray*}
&&\qquad\qquad \int_{4\eps^2}^\infty\int_{\R^3}F^\delta_\eps(t,x)\,dx\,dt\\
&&=\int_{4\eps^2}^\infty\int_{\R^3}\frac{1}{(2\eps)^5}\int_{-4\eps^2}^{0}\int_{B_{2\eps}}
F^\delta(t+s,X(t+s,t,x)+y)\,ds\,dy\,dx\,dt\\
&&=\frac{1}{(2\eps)^5}\int_{B_{2\eps}}\int_{-4\eps^2}^{0}\int_{4\eps^2}^\infty\int_{\R^3}F^\delta(t+s,X(t+s,t,x)+y)\,dx\,dt\,ds\,dy\\
&&=\frac{1}{(2\eps)^5}\int_{B_{2\eps}}\int_{-4\eps^2}^{0}\int_{4\eps^2}^\infty\int_{\R^3}F^\delta(t+s,z+y)\,dz\,dt\,ds\,dy\\
&&\leq\left(\frac{1}{(2\eps)^5}\int_{B_{2\eps}}\int_{-4\eps^2}^{0}\,ds\,dy\right)
\int_0^\infty\int_{\R^3}F^\delta(\underline{t},\underline{z})\,d\underline{z}\,d\underline{t}\\
&&=\int_0^\infty\int_{\R^3}\left(|M((-\Delta)^{-\delta/2}\nabla^2 P)|^{1+\gamma}+|\nabla u|^{2}+|\nabla^2 P|\right)\,dx\,dt.
\end{eqnarray*}
In the second equality, we have used Fubini, in the third we have
used the fact that $X$ is an incompressible flow. In the fourth
equality we did the change of variable in $(t,z)$
$$
\underline{t}=t+s\qquad \underline{z}=y+z.
$$
We then find, thanks to Tchebychev inequality,
$$
\left|\left\{F^\delta_\eps(t,x)\geq\frac{\eta^*\eps^\delta}{2(2\eps)^4}\right\}\right|\leq
2^5\frac{\int_0^\infty\int_{\R^3}F^\delta_\eps(t,x)\,dx\,dt}{\eta^*}\eps^{4-\delta}.
$$
We conclude thanks to Lemma \ref{lemmmaximalenergy}.\qquad
\qed
\vskip0.3cm

We fix $\delta>0$.  For any fixed
$(t,x)\in \Omega^\delta_\eps$ with $t\geq 4\eps^2$, we define $\ve, \Pe$,
(depending on this fixed point $(t,x)$) as  functions of two local
new variables $(s,y)\in Q_2$:
\begin{eqnarray}
\nonumber&&\ve(s,y)=\eps u(t+\eps^2 s, X(t+\eps^2 s,t,x)+\eps y)\\
\label{eq_veps}&&\qquad\qquad\qquad -\eps
u_\eps(t+\eps^2 s, X(t+\eps^2 s,t, x)),\\
 && \nonumber \Pe(s,y)=\eps^2 P(t+\eps^2 s, X(t+\eps^2 s,t,x)+\eps
y)\\
&&\label{eq_Peps}\qquad\qquad\qquad+\eps y\partial_s[ u_\eps(t+\eps^2 s, X(t+\eps^2 s,t, x))].
\end{eqnarray}
We have the following proposition.
\begin{prop}\label{prop_first estimates}
The function  $(\ve,\Pe)$ is solution to
(\ref{eq_NS}) (\ref{eq_NSenergy}) for $(s,y)\in (-4,0)\times\R^3$.
It verifies:
\begin{eqnarray}
&& \int_{\R^3}\phi(y) v_\eps(s,y)\,dy=0,\qquad s\geq -4, \label{star1V}\\
&& \int_{-4}^0\int_{B_2}|\nabla v_\eps|^2\,dy\,ds\leq \eta^*,\\
&& \int_{-4}^0\int_{B_2}|\nabla^2 P_\eps|\,dy\,ds\leq \eta^*,\\
&&\int_{-4}^0\int_{B_2}|M[(-\Delta)^{-\delta/2}\nabla^2 \Pe]|^{1+\gamma}\,dy\,ds\leq \eta^*.
\end{eqnarray}
\end{prop}
\noindent{\bf Proof.}
The fact that $(\ve,\Pe)$ is solution to
(\ref{eq_NS}) (\ref{eq_NSenergy}) and verifies (\ref{star1V}) comes from its definition (\ref{eq_veps}), (\ref{eq_Peps}), (\ref{etoile1}) and (\ref{etoile2}).
We have
\begin{equation}\label{eq_control_local}
\begin{array}{l}
\dis{\qquad\qquad\int_{Q_2}(|\nabla v_\eps|^{2}+|\nabla^2 P_\eps|)\,dy\,ds+\int_{Q_2}|M[(-\Delta)^{-\delta/2}\nabla^2 \Pe]|^{1+\gamma}\,dy\,ds}\\[0.3cm]
\dis{=\int_{Q_2}\left(\eps^4(|\nabla u|^2+|\nabla^2 P|)+\eps^{(4-\delta)(1+\gamma)}|M[(-\Delta)^{-\delta/2}\nabla^2 P]|^{1+\gamma}\right)}\\[0.3cm]
\qquad\qquad\qquad\qquad\qquad\qquad\qquad\qquad\dis{(t+\eps^2 s,X(t+\eps^2 s,t,x)+\eps y)\,dy\,ds}\\[0.3cm]
\dis{\leq\frac{1}{\eps^{1+\delta}}\int_{t-4\eps^2}^{t}\int_{B_{2\eps}}(|\nabla u|^{2}+|\nabla^2P|+M[(-\Delta)^{-\delta/2}\nabla^2 P]^{1+\gamma})}\\[0.3cm]
\qquad\qquad\qquad\qquad\qquad\qquad\qquad\qquad\qquad\qquad\qquad\dis{(s,X(s,t,x)+y)\,ds\,dy}\\[0.3cm]
\leq \eta^*.
\end{array}
\end{equation}

In the first equality, we used the definition of $\ve$ and $\Pe$, in
the second,  we used the change of variable $(t+\eps^2 s, \eps y)\to
(s,y)$ (together with the fact that $\delta<4$ and $\gamma\geq0$), and the last inequality comes from the fact that $(s,y)$ lies in $\Omega^\delta_\eps$.
\qquad\qed
\vskip0.3cm
\noindent
Our aim is to apply proposition \ref{prop_main} to $\ve$. It will be a consequence of the following section.

\section{Local study}

This section is dedicated to the following Proposition.

\begin{prop}\label{prop_estimates}
For any $\gamma>0$ and any $0<\delta<1$,
there exists a constant $\overline{\eta}<1$, and a sequence of constants $\{C_n\}$ such that
for any solution $(u,P)$ of (\ref{eq_NS}) (\ref{eq_NSenergy}) in $Q_2$
verifying
\begin{eqnarray}
&& \int_{\R^3}\phi(y) u(t,x)\,dx=0,\qquad t\geq -4, \label{Hyppropmoyenne nul+}\\
&& \int_{-4}^0\int_{B_2}|\nabla u|^2\,dx\,dt\leq \overline{\eta}, \label{eq_gradient L2}\\
&& \int_{-4}^0\int_{B_2}|\nabla^2 P|\,dx\,dt\leq \overline{\eta}, \label{eq_pression}\\
&&\int_{-4}^0\int_{B_2}|M[(-\Delta)^{-\delta/2}\nabla^2 P]|^{1+\gamma}\,dx\,dt\leq \overline{\eta}, \label{eq_local gradiant}
\end{eqnarray}
the velocity $u$ is infinitely differentiable in $x$ at $(0,0)$ and
$$
|\nabla^n u(0,0)|\leq C_n.
$$
\end{prop}

\noindent{\bf Proof.} We want to apply Proposition \ref{prop_main}. Then, by a bootstrapping argument we will get uniform controls on higher derivatives. For this, we first
 need a control of $u$ in $L^\infty(L^2)$ and a control on $P$ in $L^{\gamma+1}(L^1)$. The equation is on $\nabla P$ (not the pressure itself). Therefore, changing $P$ by $P-\int_{B_2}\phi P\,dx$ we can assume without loss of generality that
$$
\int_{\R^3}\phi(x)P(t,x)\,dx=0,\qquad -4<t<0.
$$
 To get a control in $L^{1+\gamma}(L^1)$ on the pressure it is then enough to control $\nabla P$.

\vskip0.2cm\noindent{\bf Step 1:  Control on $u$  in $L^\infty(L^{3/2})$ in $Q_{3/2}$.}
Thanks to Hypothesis (\ref{Hyppropmoyenne nul+}), there exists a constant $C$, depending only on $\phi$,  such that for any $-4<t<0$
\begin{equation}\label{eq_Poincarreu}
\|u(t)\|_{L^6(B_2)}\leq C\|\nabla u(t)\|_{L^2(B_2)}.
\end{equation}
So
$$
\|(u\cdot\nabla)u\|_{L^1(-4,0;L^{3/2}(B_2))}\leq C \|\nabla u\|^2_{L^2(Q_2)}\leq C\overline{\eta}.
$$
We need the same control on $\nabla P$.
First, multiplying (\ref{eq_NS}) by $\phi(x)$, integrating in $x$, and using  Hypothesis (\ref{Hyppropmoyenne nul+}), we find for any $-4<t<0$
\begin{equation}\label{etoile3}
\int\phi(x)(u\cdot\nabla)u\,dx+\int\phi(x)\nabla P\,dx-\int \Delta\phi u\,dx=0.
\end{equation}
So
$$
\left\|\int\phi(x)\nabla P\,dx\right\|_{L^1(-4,0)}\leq C\left(\|\nabla u\|^2_{L^2(Q_2)}+\|u\|_{L^2(-4,0;L^6(B_2))}\right)\leq C\sqrt{\overline{\eta}}.
$$
But, as for $u$,
$$
\left\|\nabla P-\int\phi\nabla P\,dx \right\|_{L^1(-4,0;L^{3/2}(B_2))}\leq
C \|\nabla^2 P\|_{L^1(Q_2)}.
$$
So, finally
\begin{equation}\label{eq_unablauL3/2}
\||(u\cdot\nabla)
u|+|\nabla P|\|_{L^1(-4,0;L^{3/2}(B_2))}\leq C\sqrt{\overline{\eta}}.
\end{equation}
Note that
\begin{eqnarray*}
&&\frac{3}{2}\frac{u }{|u |^{1/2}}\dt u =\frac{3}{2}\frac{1}{|u |^{1/2}}\dt \frac{|u |^2}{2}\\
&&\qquad\qquad
=\frac{3}{2}|u |^{1/2}\dt|u |=\dt|u |^{3/2},\\
&&\frac{3}{2}\frac{u }{|u |^{1/2}}\Delta u =\frac{3}{2}\Div\left(\frac{u }{|u |^{1/2}}\nabla u \right)
-\frac{3}{2}\frac{|\nabla u |^2}{|u |^{1/2}}+\frac{3}{4}\frac{|\nabla |u ||^2}{|u |^{1/2}}\\
&&\qquad\qquad\leq \Delta |u |^{3/2},
\end{eqnarray*}
since $|\nabla u |\geq |\nabla|u ||$.

We consider $\psi_1\in C^\infty(\R^4)$
 a nonnegative function compactly supported in $Q_{2}$ with $\psi_1=1$ in $Q_{3/2}$ and
$$
|\nabla_{t,x}\psi_1|+|\nabla_{t,x}^2\psi_1|\leq C.
$$
Multiplying (\ref{eq_NS}) by $(3/2)\psi_1(t,x)u /|u |^{1/2}$ and integrating in $x$ gives
\begin{eqnarray*}
&&\qquad\frac{d}{dt}\int\psi_1(t,x)|u |^{3/2}\,dx\\
&&\leq \int(|\dt \psi_1|+|\Delta \psi_1|)|u |^{3/2}\,dx\\
&&\qquad\qquad+\frac{3}{2}\|\psi_1^{1/3}|u |^{1/2}\|_{L^3(\R^3)}\|\psi_1^{2/3}((u \cdot\nabla)u +\nabla P)\|_{L^{3/2}(B_2)}\\
&&\leq \int(|\dt \psi_1|+|\Delta \psi_1|)|u |^{3/2}\,dx\\
&&\qquad\qquad+\frac{3}{2}\left(\int\psi_1(t,x)|u |^{3/2}\,dx\right)^{1/3}\|((u \cdot\nabla)u +\nabla P)\|_{L^{3/2}(B_2)}\\
&&\leq \alpha(t)\left(1+\int\psi_1(t,x)|u |^{3/2}\,dx\right),
\end{eqnarray*}
with
$$
\alpha(t)=\int(|\dt \psi_1|+|\Delta \psi_1|)|u |^{3/2}\,dx+\frac{3}{2}\|((u \cdot\nabla)u +\nabla P)\|_{L^{3/2}(B_2)}.
$$
Thanks to (\ref{eq_Poincarreu}) and (\ref{eq_unablauL3/2})
$$
\|\alpha\|_{L^1(-4,0)}\leq C\sqrt{\overline{\eta}}.
$$
Denoting $Y(t)=1+\int\psi_1(t,x)|u |^{3/2}\,dx$, we have
$$
\dot Y\leq \alpha Y,\qquad Y(-4)=1.
$$
Gronwall's lemma gives that for any $-4<t<0$ we have
$$
Y(t)\leq exp\left(\int_{-4}^t\alpha(s)\,ds\right).
$$
Hence, for $\overline{\eta}$ small enough:
\begin{equation}\label{eq_uL^3/2}
\|u \|_{L^{\infty}(-(3/2)^2,0;L^{3/2}(B_{3/2}))}\leq C {\overline{\eta}}^{1/3}.
\end{equation}

\vskip0.2cm\noindent{\bf Step 2:  Control on $u$  in $L^\infty(L^{2})$ in $Q_1$.}

We consider $\psi_2\in C^\infty(\R^4)$
 a nonnegative function compactly supported in $Q_{3/2}$ with $\psi_2=1$ in $Q_{1}$ and
$$
|\nabla_{t,x}\psi_2|+|\nabla_{t,x}^2\psi_2|\leq C.
$$
Multiplying inequality (\ref{eq_NSenergy}) by $\psi_2$ and integrating in $x$ gives
\begin{eqnarray*}
&&\qquad\qquad \frac{d}{dt}\left(\int\psi_2\frac{|u |^2}{2}\,dx\right)\\
&&\leq \int u \cdot\nabla\psi_2\left(\frac{|u |^2}{2}+P\right)\,dx+\int(\dt\psi_2+\Delta\psi_2)\frac{|u |^2}{2}\,dx.
\end{eqnarray*}
equalities (\ref{eq_Poincarreu}) together with (\ref{eq_unablauL3/2}) and Sobolev imbedding  gives
$$
\||u |^2+P\|_{L^1(-(3/2)^2,0;L^{3}(B_{3/2}))}\leq C{\overline{\eta}}^{1/2}.
$$
Together with (\ref{eq_uL^3/2}), this gives that
\begin{equation}\label{eq_uL^2}
\|u \|_{L^\infty(-1,0;L^2(B_1))}\leq C{\overline{\eta}}^{1/4}.
\end{equation}

\vskip0.3cm\noindent{\bf Step 3. $L^\infty$ bound in $Q_{1/2}$.} We need now to get better integrability in time on the pressure.

From (\ref{etoile3}) and (\ref{eq_uL^2}), we get
$$
\left\|\int\phi(x)\nabla P\,dx\right\|_{L^2(-1,0)}\leq C\sqrt{\overline{\eta}}.
$$
With Lemma \ref{lemm_local gradient control} and (\ref{eq_local gradiant}), this gives for $\gamma<1$
$$
\|\nabla P\|_{L^{1+\gamma}(-1,0;L^1(B_1))}\leq C\sqrt{\overline{\eta}}.
$$
Together with (\ref{eq_uL^2}), (\ref{eq_gradient L2}), and
Proposition \ref{prop_main}, this shows that
 for $\overline{\eta}$ small enough, we have
$$
|u|\leq 1\qquad\mathrm{in }\  \ Q_{1/2}.
$$

\vskip0.2cm\noindent{\bf Step 4: Obtaining more regularity.}
We now obtain higher derivative estimates by a standard bootstrapping method. We give the details carefully to ensure that the bounds obtained are universal, that is, do not depend on the actual solution $u$.

For $n\geq 1$ we define $r_n=2^{-n-3}$, $\B_n=B_{r_n}$ and $\Q_n=Q_{r_n}$.
We denote also  $\pn$ such that $0\leq\pn\leq 1$, $\pn\in C^\infty(\R^4)$,
\begin{eqnarray*}
\pn(t,x)&=&1\qquad (t,x)\in \Q_n,\\
&=&0\qquad (t,x)\in \Q_{n-1}^c.
\end{eqnarray*}

For every $n$ we have
\begin{equation}\label{eq_derive}
\dt\nabla^n u+\Div A_n+\nabla R_n-\Delta\nabla^n u=0,
\end{equation}
with
$$
A_n=\nabla^n (u\otimes u),\qquad R_n=\nabla^n P.
$$
So we have
\begin{equation}\label{eq_estimates A_n}
\|A_n\|_{L^p(\Q_{n-1})}\leq C_n\|u\|^2_{L^{2p}(-r^2_{n-1},0;W^{n,2p}(\B_{n-1}))}
\end{equation}
and thanks to Lemma \ref{lemm_decomposition of pressure}, we can split $R_n$ as
$$
R_n=R_{1,n}+R_{2,n},
$$
with
\begin{eqnarray}\label{eq_estimates R_1n}
&&\|R_{1,n}\|_{L^p(\Q_{n-1})}\leq C_n \|A_n\|_{L^p(\Q_{n-2})},\\
\nonumber
&&\|R_{2,n}\|_{L^1(-r^2_{n-1},0;W^{2,\infty}(\B_{n-1}))}\leq C_n \left(\|A_n\|_{L^p(\Q_{n-2})}+\|\nabla P\|_{L^1(\Q_{n-2})}\right)\\
\label{eq_estimates R_2n}
&&\qquad\qquad\qquad\qquad\qquad\qquad\leq C_n \left(\|A_n\|_{L^p(\Q_{n-2})}+1\right).
\end{eqnarray}

Moreover we have:
\begin{eqnarray*}
&&\dt(\pn \nabla^n u)-\Delta(\pn \nabla^n u)\\
&&\qquad=-\Div(A_n \pn)+\nabla \pn A_n\\
&&\qquad\qquad -\nabla(\pn R_n)+(\nabla\pn) R_n\\
&&\qquad\qquad +\Delta\pn \nabla^n u-2\Div(\nabla\pn\nabla^n u)\\
&&\qquad\qquad +(\dt\pn)\nabla^n u.
\end{eqnarray*}

Note that $\pn\nabla^n u=0$ on $\partial \Q_{n-1}$.
So
\begin{equation}\label{eq_v12}
\pn\nabla^n u=V_{1,n}+V_{2,n}
\end{equation}
with
\begin{eqnarray*}
&&\dt V_{1,n}-\Delta V_{1,n}=-\Div(A_n \pn)+\nabla \pn A_n\\
&&\qquad -\nabla(\pn R_{1,n})+(\nabla\pn) R_{1,n}\\
&&\qquad +\Delta\pn \nabla^n u-2\Div(\nabla\pn\nabla^n u)\\
&&\qquad +(\dt\pn)\nabla^n u\\
&&\qquad\qquad\qquad=F_n,\\
&&V_{1,n}=0\qquad \mathrm{for } \ t=-r_{n-1}^2,
\end{eqnarray*}
and
\begin{eqnarray*}
&&\dt V_{2,n}-\Delta V_{2,n}=-\nabla(\pn R_{2,n})+R_{2,n}(\nabla \pn),\\
&&V_{2,n}=0\qquad \mathrm{for } \ t= -r_{n-1}^2.
\end{eqnarray*}
Thanks to (\ref{eq_estimates A_n}) and (\ref{eq_estimates R_1n}), we have
$$
\|F_n\|_{L^{p}(-r^2_{n-1},0;W^{-1,p}(\B_{n-1}))}\leq C_n\left(1+\|u\|^2_{L^{2p}(-r^2_{n-2},0;W^{n,2p}(\B_{n-2}))}\right).
$$
So, from Lemma \ref{lemm_parabolic regularity},
\begin{eqnarray*}
&&\|V_{1,n}\|_{L^{p}(-r^2_{n-1},0;W^{1,p}(\B_{n-1}))}\leq C\|F_n\|_{L^{p}(-r^2_{n-1},0;W^{-1,p}(\R^3))},\\
&&\|V_{2,n}\|_{L^{\infty}(-r^2_{n-1},0;W^{1,\infty}(\B_{n-1}))}\leq C\|\pn\nabla R_{2,n}\|_{L^{1}(-r^2_{n-1};W^{1,\infty}(\R^3))}\\
&&\qquad\qquad\qquad\qquad+C\| R_{2,n}(\nabla \pn)\|_{L^{1}(-r^2_{n-1}W^{1,\infty}(\R^3))}\\
&&\qquad\qquad\leq C_n \left(1+\|u\|^2_{L^{2p}(-r^2_{n-2},0;W^{n,2p}(\B_{n-2}))}\right),
\end{eqnarray*}
where we have used (\ref{eq_estimates A_n}) and (\ref{eq_estimates R_2n})  in the last line.

Hence, from (\ref{eq_v12}) and using that $\pn=1$ on $\Q_n$, we have for any $1<p<\infty$
$$
\|\nabla^n u\|_{L^p(-r^2_n,0;W^{1,p}(\B_n))}\leq C_n \left(1+\|u\|^2_{L^{2p}(-r^2_{n-2},0;W^{n,2p}(\B_{n-2}))}\right).
$$
By induction we find that for any $n\geq 1$, and any $1\leq p<\infty$, there exists a constant $C_{n,p}$ such that
$$
\|u\|_{L^{2^{-n}p}(-r^2_n,0;W^{n,2^{-n}p}(\B_n))}\leq C_{n,p}.
$$
This is true for any $p$, so for $n$ fixed, taking $p$ big enough and using
 Sobolev imbedding, we show that  for any $1\leq q<\infty$, there exists a constant $C_{n,q}$ such that
$$
\|u\|_{L^q(-r^2_{n+1},0; W^{n,\infty}(\B_{n+1}))}\leq C_{n,q}.
$$
As (\ref{eq_estimates A_n}), we get
that
$$
\|A_n\|_{L^1(-r^2_{n+3},0;W^{2,\infty}(\B_{n+3}))}\leq C_{n}.
$$
Thanks to Lemma \ref{lemm_decomposition of pressure}, we get
\begin{eqnarray*}
&&\|R_{1,n}\|_{L^1(-r^2_{n+4},0;W^{1,\infty}(\B_{n+4}))}\leq C_{n},\\
&&\|R_{2,n}\|_{L^1(-r^2_{n+4},0;W^{1,\infty}(\B_{n+4}))}\leq C_{n}.
\end{eqnarray*}
Hence
$$
\|\dt\nabla^n u\|_{L^1(-r^2_{n+4},0;L^{\infty}(\B_{n+4}))}\leq C_n,
$$
and finally
$$
\|\nabla^n u\|_{L^\infty(\Q_{n+4})}\leq C_n.
$$
\qquad \qed

\section{From local to global}

Let us fix $\delta>0$. We take $\eta^*\leq\overline{\eta}$ and consider any $\eps>0$ such that $4\eps^2\leq t_0$. Then from Proposition \ref{prop_estimates} and Proposition \ref{prop_first estimates}, for any $(t,x)\in \Omega^\delta_\eps\cap\{t\geq t_0\}$,
we have
$$
|\nabla^n_y v_\eps(0,0)|\leq C_n,
$$
where $v_\eps$ is defined by (\ref{eq_veps}).
But for any $n\geq1$, we have
$$
\nabla^n_y v_\eps(0,0)=\eps^{n+1}\nabla^nu(t,x).
$$
Hence
$$
\left|\left\{(t,x)\in\Omega \setminus |\nabla^nu(t,x)|\geq \frac{C_n}{\eps^{n+1}}\right\}\right|\leq |[\Omega^\delta_\eps]^c|.
$$
And thanks to Lemma \ref{lemm_main},
This measure is smaller than
$$
\frac{C}{\eta^*}\left(\|u^0\|^2_{L^2(\R^3)}+\|u^0\|_{L^2(\R^3)}^{2(\gamma+1)}\right)\eps^{4-\delta}.
$$
We denote
$$
R=\left(1+\frac{4}{t_0}\right)^{\frac{n+1}{2}}.
$$
For $k\geq1$, we use our estimate with $\eps^{n+1}=R^{-k}$ to get
$$
\left|\left\{(t,x)\in\Omega \setminus \frac{|\nabla^nu(t,x)|}{C_n}\geq R^k\right\}\right|\leq \frac{C\left(1+\|u^0\|_{L^2(\R^3)}^{2(\gamma+1)}\right)}{R^{k\frac{4-\delta}{n+1}}}.
$$
So, for $p<\frac{4-\delta}{n+1}$
\begin{eqnarray*}
&&\left\|\frac{\nabla^n u}{C_n}\right\|^p_{L^p(\Omega)}\leq \left|\left\{(t,x)\in\Omega \setminus \frac{|\nabla^nu(t,x)|}{C_n}\leq R\right\}\right| R^p\\
&&\qquad\qquad +\sum_{k=1}^\infty R^{(k+1)p}\left|\left\{(t,x)\in\Omega \setminus \frac{|\nabla^nu(t,x)|}{C_n}\geq R^k\right\}\right|\\
&&\qquad \leq |\Omega|R^p+CR^p\left(1+\|u^0\|_{L^2(\R^3)}^{2(\gamma+1)}\right)\sum_{k=1}^\infty R^{k\left(p-\frac{4-\delta}{n+1}\right)}\\
&&\leq |\Omega|R^p+\frac{CR^p}{1-R^{p-\frac{4-\delta}{n+1}}}\left(1+\|u^0\|_{L^2(\R^3)}^{2(\gamma+1)}\right).
\end{eqnarray*}
The results holds for any $\delta>0$ which ends the proof of Theorem \ref{theo_main}.\qquad \qed

\vskip0.2cm\noindent{\bf Acknowledgment}: This work was partially supported by NSF Grant DMS-0607053. We thank Prof. Caffarelli for many insightful discussions and advices.

\bibliography{biblio}

\end{document}